\documentclass[11pt,a4paper]{amsart}
\usepackage{amsfonts}
\usepackage{}
\usepackage{amsfonts}
\usepackage{amsmath,xypic}
\usepackage{mathrsfs}
\usepackage{amssymb}

\usepackage{CJK,CJKnumb}
\usepackage[CJKbookmarks,colorlinks,
            linkcolor=black,
            anchorcolor=black,
            citecolor=blue]{hyperref}
\usepackage{color}              % 支持彩色
\usepackage{indentfirst}        %首行缩进宏包
\usepackage{latexsym,bm}        % 处理数学公式中和黑斜体的宏包
\usepackage{amsmath,amssymb}    % AMSLaTeX宏包 用来排出更加漂亮的公式
\usepackage{pstricks}
\usepackage{pst-node}
\usepackage{pst-tree}
\usepackage{pst-plot}
\usepackage{pst-text}
\usepackage{graphicx}
\usepackage{cases}
\usepackage{pifont}
\usepackage{txfonts}
\usepackage[all,knot,poly]{xy}

\allowdisplaybreaks[4]  % \eqnarray如果很长，影响分栏、换行和分页
\CJKtilde   %用于解决英文字母和汉字的间距问题。例如：变量~$x$~的值。

\newcommand{\al}{\alpha}
\newcommand{\be}{\beta}

\newcommand{\De}{\Delta}

\newcommand{\mg}{\mathfrak{g}}

\newcommand{\bK}{\mathbb{K}}
\newcommand{\bZ}{\mathbb{Z}}

\newcommand{\ot}{\otimes}

\newcommand{\si}{\sigma}

\newcommand{\ve}{\varepsilon}

\newcommand{\lan}{\langle}
\newcommand{\ran}{\rangle}

\newcommand {\fr}{\frac}
\newcommand{\A}{\mathcal{A}}

\newcommand{\W}{\mathcal{W}}
\newcommand{\Z}{\mathbb{Z}}
\newcommand{\p}{\partial}
\newcommand{\bi}{\beta_i}
\newcommand{\vei}{\varepsilon_i}
\newcommand{\T}{\Theta}

\begin{document}
\numberwithin{equation}{section}

\newtheorem{theorem}{Theorem}[section]

\newtheorem{lem}[theorem]{Lemma}

\newtheorem{cor}[theorem]{Corollary}
\newtheorem{prop}[theorem]{Proposition}

\theoremstyle{remark}
\newtheorem{rem}[theorem]{Remark}

\newtheorem{defn}[theorem]{Definition}

\newtheorem{exam}[theorem]{Example}

%%%%%%%%%%%%%%%%%%%%%%%%%%%%%%%%%%%%%%%%%%%%%%%
%%%%%%%%%%%%%%%%%%%%%%%%%%%%%%%%%%%%%%%%%%%%%%%%%
\title[Matching realization of $U_q(sl_{n+1})$ in $\mathcal W_q(2n)$]
{Matching realization of $U_q(sl_{n+1})$ of higher rank in the quantum Weyl algebra $\mathcal W_q(2n)$}

\author[Hu]{Naihong Hu}
\address{Department of Mathematics,
Shanghai Key Laboratory of Pure Mathematics and Mathematical
Practice, East China Normal University, Minhang Campus, Dong Chuan
Road 500, Shanghai 200241, PR China} \email{nhhu@math.ecnu.edu.cn}
\thanks{The programme is supported by the NNSF of China (No. 11271131).}

\author[Wang]{Shenyou Wang}
\address{Department of Mathematics, Shanghai Key Laboratory of Pure Mathematics and Mathematical Practice,
East China Normal University,
Minhang Campus, Dong Chuan Road 500, Shanghai 200241, PR China}
\email{youxiong198@163.com}

\subjclass{Primary 17B37; Secondary  81R50}

\date\today
\begin{abstract}
In the paper, we further realize the higher rank quantized universal
enveloping algebra $U_q(sl_{n+1})$ as certain quantum differential
operators in $\W_q(2n)$ defined over the quantum divided power
algebra $\mathcal{A}_q(n)$ of rank $n$. We give the quantum
differential operators realization for both the simple root vectors
and the non-simple root vectors of $U_q(sl_{n+1})$. The nice
behavior of the quantum root vectors formulas under the action of
the Lusztig symmetries once again indicates that our realization
model is naturally matched.
\end{abstract}

\maketitle
\section{Introduction and Preliminaries}
\subsection{}
It is well known that there exist quantum analogues $U_q(\mg)$
 for Lie algebras $\mg$ with (generalized) Cartan matrices.
 But how to quantize the Lie algebras of Cartan type such that their $0$-component parts coincide with the standard Drinfeld-Jimbo quantum groups is still an
 open problem. This is a hard question aimed by the first author when he was a Humbolt research fellow almost twenty years ago. On the other hand, fortunately,
 a series work on the mudular quantizations of the modular simple restricted Lie algebras of Cartan type have been done in recent years via the Jordanian twists and modular reductions, and these provide us with some new
 non-pointed Hopf algebras of prime-power dimensions defined over a field of positive characteristic, see \cite{HW1}, \cite{HW2}, \cite{HTW}. This will be significant if one considers the Kaplansky's 10 questions
 proposed in early 1975 which are related to classifying some finite-dimensional Hopf algebras in some sense, and compares with the seminal work on classifying finite-dimensional complex pointed Hopf algebras
 due to Andruskiewitsch-Schneider \cite{AS} in 2010. Here we would like to point out that these modular quantizations are non-pointed and do not belong to
 the standard ones in the sense of Drinfeld-Jimbo type.

 As the Cartan type Lie algebras contain four series $W$, $S$, $H$,
 $K$, our concern will focus on the Jacobson-Witt type for trying to attacking this problem as a first
 step in this paper.

\subsection{}
Let
$\A(n)=\text{Span}_\mathbb{F}\{x^{(\al)}\mid\al\in\mathbb{Z}^n_+\}$
 be the divided power algebra. Its multiplication is defined by
$$x^{(\al)}x^{(\be)}=\binom {\al+\be}{\al}x^{(\al+\be)},$$
where $\binom {\al+\be}{\al}=\prod\limits_{i=1}^n\binom
{\al_i+\be_i}{\al_i},~ \binom
{\al_i+\be_i}{\al_i}=\frac{(\al_i+\be_i)!}{\al_i!\be_i!}$.

The Jacobson-Witt algebra $\W(n)$ is defined to be
 the derivation algebra of divided power algebra $\A(n)$ (see~\cite{SF} for more details):
\begin{eqnarray*}
\W (n)=Der\A(n)=\left.\left\{\,\sum\limits_{i=1}^n
f_i\p_i\,\right|\,f_i\in\A(n)\, \right\}.
\end{eqnarray*}

%An open question is how to quantize~$\W(n)$, namely,
%how to find a quantum algebra~$U_q(\W(n))$
%such that it is exactly the quantum analogue of~$\W(n)$ ?
%~(see~\cite{Hu2}~for the case $n=1$)

Set $deg(x^{(\al)})=|\al|=\sum_{i=1}^n\al_i,deg(\p_i)=-1$.
Then there exists a natural gradation over $\W(n)$:
$$\W(n)=\bigoplus_{i\geq -1} \W(n)_i=\bigoplus_{i\geq -1} \text{Span}_\mathbb{F}\left.\left\{\,x^{(\al)}\p_i\,\right|\,|\al|=i+1\,\right\},$$
where $\W(n)_i$ is the subspace of the derivations of degree $i$.

Note that $\W(n)_0\cong gl_n=sl_n\oplus \mathbb{C}c$,
where $c=diag(1,\cdots,1)$
and the isomorphism is obtained by mapping the differential operators $x^{(\vei)}\p_j$'s to the matrix units $E_{ij}$'s, where the action of the operator $x^{(\vei)}\p_j$
acting on $x^{(\al)}$ is given by $x^{(\vei)}\p_j(x^{(\al)})=\binom{\al{+}\vei{-}\varepsilon_j}{\vei}x^{(\al{+}\vei{-}\varepsilon_j)}$.
This means that $gl_n$ is realized
as certain differential operators defined over the divided power algebra $\A(n)$.

\subsection{}
It is natural to ask how to realize $U_q(gl_n)$ using quantum
differential operators such that $U_q(\W(n)_0)\cong U_q(gl_n)$? To
deal with this problem, the first author \cite{Hu1} introduced the
quantum version of $\A(n)$, namely, the quantum divided power algebra
$\A_q(n)$. As a vector space, it is generated by
$\{\,x^{(\alpha)}\mid\alpha\in\Z^n_+\,\}$ and its multiplication is
given by:
\begin{eqnarray*}
x^{(\alpha)}x^{(\beta)}=q^{\alpha*\beta}
\begin{bmatrix}
\alpha+\beta\\ \alpha
\end{bmatrix}
x^{(\alpha+\beta)},
\end{eqnarray*}
where
$\begin{bmatrix}
\alpha{+}\beta\\ \alpha
\end{bmatrix}=\prod\limits_{i=1}^n\begin{bmatrix}
\alpha_i{+}\beta_i\\ \alpha_i\end{bmatrix},~~
\begin{bmatrix}
\alpha_i{+}\beta_i\\
\alpha_i\end{bmatrix}=\fr{[\alpha_i{+}\beta_i]!}{[\alpha_i]![\beta_i]!},~~
[\alpha_i]!=[\al_i][\al_i{-}1]\cdots [1]$,
$[\al_i]=\fr{q^{\al_i}-q^{-\al_i}}{q-q^{-1}}$, $(\alpha_i,
\beta_i\in\Z_+)$ and for the definition of $\alpha*\beta$,
see \cite{Hu1} (also see subsection 1.7). We will write $x^{(\vei)}$ simply $x_i$ when no
confusion arises.

Also, the first author \cite{Hu1} introduced the quantum differential operators and
algebra automorphisms defined over $\A_q(n)$ as follows.

For $1\le i \le n$, $\fr{\p_q}{\p x_i}$ is defined as the special
$q$-derivatives over $\A_q(n)$ by
\begin{eqnarray}\label{partial}
\fr{\p_q}{\p x_i}(x^{(\beta)})=q^{-\ve_i*\beta}(x^{(\beta-\ve_i)}),\quad
\forall x^{(\beta)}\in\A_q(n).
\end{eqnarray}
For simplicity of notation, we write $\p_i$ for $\fr{\p_q}{\p {x_i}}$.

For $1\le i \le n$ and $\alpha \in \Z^n_+$,
 define the algebra automorphisms $\si_i$
and $\T(\alpha)$ of quantum divided power algebra as follows:
\begin{eqnarray}
&&\si_i (x^{(\beta)})=q^{\bi}(x^{(\beta)}),\quad \forall
x^{(\beta)}\in\A_q,\\\label{sigma}
&&\T(\alpha)(x^{(\beta)})=\theta(\alpha,\beta)(x^{(\beta)}),\quad
\forall x^{(\beta)}\in\A_q.\label{theta}
\end{eqnarray}
Taking $q=1$, we have $\si_i=\text{id}=\T(\alpha).$

These quantum differential operators and automorphisms yield a Hopf
algebra $D_q$ (see subsection 1.8 for more details) such that the
quantum divided power algebra  $\A_q$ is a $D_q$-module algebra in
the sense of \cite{Sw}. This allows us to make their smash product
algebra $\A_q\sharp D_q$ and define it to be the quantum Weyl
algebra (see \cite{Hu1}), denoted by $\mathcal{W}_{q}(2n)$, which is
different from those that have appeared in the literature (for
instance \cite{DF}, \cite{Ha}, \cite{Ma}, etc.).

Although $\mathcal{W}_{q}(2n)$ itself is not a Hopf algebra, it
contains many interesting Hopf subalgebras such as $U_q(gl_n)$ and $U_q(sl_n)$
(these can be realized as some quantum differential operators in
$W_q(2n)$), and also provides some strong smash product algebras in rank $1$ case (revised version for root of unity)
for calculating their cyclic homologies (see \cite{ZH}).

%In \cite{Hu1}, $U_{q}(sl_n)$ is realized as certain~$q$-differential
%operators in $\W_q(2n)$
%defined over the quantum divided algebra $\A_q(n)$.
\begin{prop}$($\cite[Thm. 4.1]{Hu1}$)$\label{qsln}
 For any monomial $x^{(\beta)}\in \mathcal{A}_q(n)$ and $1\le i\le n-1$,
 set
\begin{eqnarray*}
e_i(x^{(\beta )})&=&x_i{\partial}_{i+1}{\sigma}_i(x^{(\beta)}),\\
f_i(x^{(\beta)})&=&{{\sigma}_i}^{-1}
x_{i+1}{\partial}_{i}(x^{(\beta)}),\\
K_i(x^{(\beta)})&=&{\sigma}_i{\sigma}_{i+1}^{-1}(x^{(\beta)}),\\
K_i^{-1}(x^{(\beta)})&=&{\sigma}_i^{-1}{\sigma}_{i+1}(x^{(\beta)}).
\end{eqnarray*}
These formulas define the structure of a
$U_{q}(sl_n)$-module algebra on  $\A_q(n)$.
\end{prop}

%The following corollary gives the corresponding realization of
%$U_{q}(gl_n)$.
\begin{cor}$($\cite[Coro. 4.1]{Hu1}$)$\label{qgl_n}
 For any monomial $x^{(\beta)}\in \mathcal{A}_q(n)$, set
\begin{eqnarray*}
e_i(x^{(\beta )})&=&x_i{\partial}_{i+1}{\sigma}_i(x^{(\beta)})\quad(1\le i \le n-1),\\
f_i(x^{(\beta)})&=&{{\sigma}_i}^{-1}
x_{i+1}{\partial}_{i}(x^{(\beta)})\quad(1\le i \le n-1),\\
k_i(x^{(\beta)})&=&{\sigma}_i(x^{(\beta)})\quad(1\le i \le n),\\
k_i^{-1}(x^{(\beta)})&=&{\sigma}_i^{-1}(x^{(\beta)})\quad(1\le i \le n).
\end{eqnarray*}
These formulas define the structure of a $U_{q}(gl_n)$-module
algebra over $\A_q(n)$.
\end{cor}
Actually, this realizes the $U_q(\W(n)_0)\cong U_q(gl_n)$ as some
quantum differential operators in $n$ variables defined over
$\A_q(n)$.

\subsection{}
Let
$L_1=\text{Span}_\mathbb{F}\{x_k\sum\limits_{i=1}^nx_i\p_i\mid
1\le k \le n\} \subset \W(n)_1$ be an abelian Lie subalgebra of
$\W(n)$. Set $L_{n+1}:=\W(n)_{-1}\oplus \W(n)_0\oplus L_1 \subset
\W(n)$, It is easy to check that $L_{n+1}$ is also a Lie subalgebra
of $\W(n)$, which is isomorphic to $sl_{n+1}$ by Zhao-Xu \cite{ZX}.

\begin{lem}$($\cite[Lemma 3.1.1]{ZX}$)$ \label{sln+1}
The special linear Lie algebra $sl_{n+1}$ is isomorphism to
$L_{n+1}$ with the following identification of Chevalley generators
\begin{eqnarray*}
e_i&=&x_{i}\p_{i+1},~ (1\le i \le n{-}1),\qquad\qquad\qquad
e_n=x_n\sum\limits_{i=1}^nx_i\p_i,
\\
f_i&=&x_{i{+}1}\p_{x_{i}},~ (1\le i \le n{-}1),\qquad\qquad\qquad
f_n=-\p_n,\\
h_i&=&x_i\p_i-x_{i+1}\p_{i+1},~ (1\le i \le n{-}1),\qquad\,
h_n=\sum\limits_{i=1}^nx_i\p_{x_i}+x_n\p_n.
\end{eqnarray*}
\end{lem}

This means that one can give the differential operator realization
for the higher rank special linear Lie algebra $sl_{n+1}$ using the
differential operators in $n$ variables (rather than $n{+}1$
variables).

In this paper, we will generalize the above result
 to the quantum case, that is, give the quantum differential operators
realization for the higher rank quantized enveloping algebra
$U_q(sl_{n+1})$ using certain quantum differential operators in
$\W_q(2n)$, namely, we realize the quantum algebra of Lie subalgebra
$L_{n+1}:=\W(n)_{-1}\oplus \W(n)_0\oplus L_1$ of $\W(n)$ as certain
quantum differential operators defined over quantum divided power
algebra $\A_q(n)$. Hence we go ahead towards quantizing the entire
Lie alegebra  $\W(n)$. We emphasize that we are able to  realize the
higher rank quantum group $U_q(sl_{n+1})$ using  quantum
differential operators in $n$ variables defined over $\A_q(n)$
rather than $n{+}1$ variables defined over $\A_q(n{+}1)$, since the
latter follows directly from Proposition \ref{qsln}.

%There are two main results in this paper.
%The one mentioned above is to give a kind of quantized realizations of
%the simple root vectors of higher rank $U_q(sl_{n+1})$;
%The other one is realizing further the non-simple root vectors
%as certain quantum differential operators in $n$ variables.

\subsection{}
In \cite{Hu1}, after giving the realization of simple root vectors
of $U_q(gl_n)$, the first author gave further the realization of all
root vectors of $U_q(gl_n)$, which is proved to be coincident with
one of four kinds of root vectors in the sense of Lusztig using the
braid automorphisms of $U_q(sl_n)$ (also known as Lusztig symmetries
\cite{Lus}).

For $1\le i <j \le n$, let~$e_{ij}=x_i\p_j\si_i$, and for $1\le j
<i\le n$, set~$e_{ij}=\si_j^{-1}x_i\p_j$.
\begin{prop}$($\cite[Prop. 4.6]{Hu1}$)$\label{prop4.6}
Identifying the generators of $U_q(sl_{n})$
with certain quantum differential operators in $\W_q(2n)$,
i.e.,
$e_i:=e_{i,i+1},f_i:=e_{i+1,i},K_i:=\sigma_i\sigma_{i+1}^{-1}$
$(1\le i\le n{-}1)$,  we have

(1) $e_{\alpha_{ij}}$ correspond to $e_{i,j}~(1\leq i<j\leq n)$,
where $e_{\alpha_{ij}}$ are the positive root vectors associated to
those positive roots~$\alpha_{ij}=\varepsilon_i-\varepsilon_j$.

(2) $f_{\alpha_{ji}}$ correspond to $e_{i,j}~(1\leq j<i\leq n)$,
where $f_{\alpha_{ji}}$ are the negative root vectors associated to
those negative roots $\alpha_{ji}=\varepsilon_i-\varepsilon_j$.
\end{prop}

 In this paper, after giving the realization of simple root vectors
 of higher rank quantum group $U_q(sl_{n+1})$ ,
 we also realize further the non-simple root vectors using the Lusztig
 symmetries and show that one of the four kinds of
 Lusztig root vectors can be specified under the realization.

% The good performance of ~Lusztig~symmetries in our realization model further
%illustrates that our realization model is natural.

\subsection{}
For the convenience of the reader,
 we summarize the relevant results without proofs to make our exposition self-contained.

The quantum group $U_q(sl_n)$ \cite{Jan} is an associative algebra over $\bK$ generated by the elements
$e_i,f_i,K_i^{\pm1}~~(1\le i\le n-1)$, subject to the relations:
\begin{eqnarray*}
(R1)&& K_iK_i^{-1}=K_i^{-1}K_i=1,~~~~~~~~~~~K_iK_j=K_jK_i,\\
 (R2)&&  K_ie_j{K_i}^{-1}=q^{a_{i,j}}e_j,~~~~~~~K_if_j{K_i}^{-1}=q^{-a_{i,j}}f_j,\\
(R3)&&\qquad \quad [e_i,f_j]=\delta_{ij}\frac{K_i-K_i^{-1}}{q-q^{-1}},\\
(R4)&&e_i^2e_j-(q+q^{-1})e_ie_je_i+e_je_i^2=0\quad(|i-j|=1),\\
(R5)&&\qquad e_ie_j=e_je_i\quad(|i-j|>1),\\
(R6)&&f_i^2f_j-(q+q^{-1})f_if_jf_i+f_jf_i^2=0\quad(|i-j|=1),\\
(R7)&&\qquad f_if_j=f_jf_i\quad(|i-j|>1).
\end{eqnarray*}
where $q\in \bK^*$ and $(a_{ij})$ is the Cartan matrix of type
$A_{n-1}$.

 $U_q(sl_n)$ has a Hopf algebra structure
with the comultiplication, the counit and the antipode given by:
\begin{eqnarray*}
&& \Delta(K_i^{\pm1})=K_i^{\pm1}\ot K_i^{\pm1},\quad
 \Delta(e_i)=e_i\ot K_i+1\ot e_i,\\
&&\Delta(f_i)=f_i\ot 1+ K_i^{-1}\ot f_i, \qquad \ve(K_i^{\pm1})=1,\\
&& \ve(e_i)=\ve(f_i)=0, \qquad \qquad S(K_i^{\pm1})=K_i^{\mp1},\\
&&S(e_i)=-e_iK_i^{-1}, \qquad\qquad S(f_i)=-K_if_i.
\end{eqnarray*}

Let $P=\text{Span}_{\bZ}\{\ve_1,\cdots,\ve_n\}$ be the weight lattice for $gl_n$
 and $\{\alpha_i=\ve_i-\ve_{i+1}\mid 1\le i \le n{-}1\}$
be the set of simple roots of $gl_n$. Define a symmetric bilinear form~
$\lan\cdot,\cdot\ran:P\times P \rightarrow \bZ$ such that
$$\lan\ve_i,\ve_j\ran=\delta_{ij}.$$
 Now we can give the presentation of~$U_q(gl_n)$ as follows:
\begin{eqnarray*}
k_ik_i^{-1}=k_i^{-1}k_i&=&1,~~~~~~~~~~~~~~~~k_ik_j=k_jk_i \quad (1\le i\le n),\\
K_i&=&k_ik_{i+1}^{-1}  \quad (1\le i \le n-1),\\
k_ie_j{k_i}^{-1}&=&q^{\lan\vei,\alpha_j\ran}e_j,
~~~~~~k_if_j{k_i}^{-1}=q^{-\lan\vei,\alpha_j\ran}f_j.
\end{eqnarray*}
and keep the remaining relations in the definition of  $U_q(sl_n)$ invariant.

\subsection{}
 We introduce a product $*$ and a mapping $\theta$
 of two integer $n$-tuples following \cite{Hu1}.
For any
~$\alpha=(\alpha_1,\cdots,\alpha_n),~\beta=(\beta_1,\cdots,\beta_n)\in\Z^n
$, define
$\alpha*\beta:=\sum\limits_{j=1}^{n-1}\sum\limits_{i>j}\alpha_i\beta_j
$ and define $\theta: \Z^n \times \Z^n\rightarrow k^{*}$ by
$\theta(\alpha,\beta)=q^{\alpha*\beta-\beta*\alpha}$.

\begin{lem}$($\cite[\S2.1]{Hu1}$)$\label{*product and theta}
 We list some useful properties of the product $*$ and the mapping $\theta$:
\begin{eqnarray*}
(\alpha+\beta)*\gamma&=&\alpha*\gamma+\beta*\gamma,\quad
\alpha*(\beta+\gamma)=\alpha*\beta+\alpha*\gamma,\\
\theta(\alpha{+}\beta,\gamma)&=&\theta(\alpha,\gamma)\theta(\beta,\gamma),\quad
\theta(\alpha,\beta{+}\gamma)=\theta(\alpha,\beta)\theta(\alpha,\gamma),\\
\theta(\alpha,0)&=&1=\theta(0,\alpha),\quad
\theta(\alpha,\beta)\theta(\beta,\alpha)=1=\theta(\alpha,\alpha).
\end{eqnarray*}
Taking $\al$ or $\be$ to be $\ve_i$ or $\ve_j$, we have
\begin{eqnarray*}
\vei*\beta=\sum\limits_{s<i}\beta_s, \quad
\beta*\vei=\sum\limits_{s>i}\beta_s \quad (1\le i\le n),\quad
\theta(\vei,\ve_j)&=&\left\{
\begin{array}{cc}
q&\quad(i>j)\\
1&\quad(i=j)\\
q^{-1}&\quad(i<j)
\end{array}
\right.
\end{eqnarray*}
where $\ve_i =(\delta_{1i},\cdots,\delta_{ni}), ~(1\le i< n)$.
\end{lem}

\subsection{}
In this section, we will give the definition of quantum Weyl algebra
$\mathcal{W}_{q}(2n)$ which is defined in \cite{Hu1} as a kind of
smash product algebras.

Let us first recall some elementary properties of the quantum
differential operators and algebra automorphisms in subsection
1.3, which will be needed in the construction of
$\mathcal{W}_{q}(2n)$.
\begin{prop}$($\cite[Prop.3.1]{Hu1}$)$\label{deau}
 (1) ~~$\Theta(\alpha)\Theta(\beta)=\Theta(\alpha+\beta),$ in particular,
$\Theta(-\alpha_{i})=\sigma_{i}\sigma_{i+1},$
where~$\alpha_{i}$ is a simple root in a root system of type $A_{n-1}$.\\
(2) ~~$\p_i$ is a $(\T(-\ve_i)\si_i^{\pm1},\si_i^{\mp1})$ derivation
of $\A_q$, that is to say,
$$\partial_{i}(x^{(\beta)}x^{(\gamma)})=\partial_{i}(x^{(\beta)})\si_i^{\mp1}(x^{(\gamma)})
+\T(-\ve_i)\si_i^{\pm1}(x^{(\beta)})\partial_{i}(x^{(\gamma)}).$$
(3)~~~~~~~~~~~~$\p_i\p_j=\theta(\ve_i,\ve_j)\p_j\p_i.$\\
(4)~~~~$x^{(\alpha)}(x^{(\beta)}x^{(\gamma)})=(x^{(\alpha)}x^{(\beta)})x^{(\gamma)}
=\theta(\alpha,\beta)x^{(\beta)}(x^{(\alpha)}x^{(\gamma)})=\T(\alpha)(x^{(\beta)})(x^{(\alpha)}x^{(\gamma)}).$\\
(5)~~~~~~~~~~~~$\si_i(x^{(\beta)}x^{(\gamma)})=\si_i(x^{(\beta)})\si_i(x^{(\gamma)}),~~~~~\T(\alpha)(x^{(\beta)}x^{(\gamma)})=\T(\alpha)(x^{(\beta)})\T(\alpha)(x^{(\gamma)}).$\\
(6) ~~$x^{(\alpha)}\p_i$ is a
$(\T(\alpha-\ve_i)\si_i^{\pm1},\si_i^{\mp1})$-derivation of
$\A_q(n)$, namely,
$$x^{(\alpha)}\partial_{i}(x^{(\beta)}x^{(\gamma)})
=x^{(\alpha)}\partial_{i}(x^{(\beta)})\si_i^{\mp1}(x^{(\gamma)})
+\T(\alpha-\ve_i)\si_i^{\pm1}(x^{(\beta)})x^{(\alpha)}\partial_{i}
(x^{(\gamma)}).$$
\end{prop}

Let $D_q$ be the associative algebra over $\bK$ generated by the
elements $\Theta(\pm\vei),\si_i^{\pm1},\p_i,(1\le i\le n)$,
satisfying the following relations:
\begin{eqnarray*}
&&\T(\pm\ve_i)\T(\mp\ve_i)=1=\si_i^{\pm1}\si_i^{\mp1},\quad
\Theta({-}\ve_i{+}\ve_{i+1})=\si_i \si_{i+1},\\
&&\Theta(\vei)\Theta(\ve_j)=\Theta(\vei{+}\ve_j)=\Theta(\ve_j)\Theta(\vei),\quad
\si_i\si_j=\si_j\si_i,\\
&&\si_i\Theta(\ve_j)=\Theta(\ve_j)\si_i,\quad
 \Theta(\vei)\p_{j}\Theta({-}\vei)=\theta(\ve_j,\vei)\p_{j},\\
&&\si_i\p_{j}\si_i^{-1}=q^{-\delta_{ij}}\p_{j},\quad
\p_i\p_j=\theta(\vei,\ve_j)\p_j\p_i.
\end{eqnarray*}
Furthermore, $D_q$ can be equipped with two~Hopf~algebra structures
with the structure mappings
$\De^{(\pm)},\ve$ and $S^{(\pm)}$ are given by
\begin{eqnarray*}
\De^{(\pm)}(\si_i^{\pm 1})&=&\si_i^{\pm 1}\ot \si_i^{\pm 1},\quad
\De^{(\pm)}(\T(\pm\ve_i))=\T(\pm\ve_i)\ot \T(\pm\ve_i),\\
\De^{(\pm)}(\p_i)&=&\p_i\ot \si_i^{\mp 1}+\T(-\ve_i)\si_i^{\pm 1}\ot \p_i,\quad
\ve(\p_i)=0.\\
\ve(\si_i^{\pm 1})&=&1=\ve(\T(\pm\ve_i)),\quad
S^{(\pm)}(\si_i^{\pm 1})=\si_i^{\mp 1},\\
S^{(\pm)}(\T(\pm \ve_i))&=&\T(\mp\ve_i),\quad
S^{(\pm)}(\p_i)=-q^{\pm 1}\T(\ve_i)\p_i.
\end{eqnarray*}

Let $D_q^{(\pm)}:=(D_q,\De^{(\pm)},\ve,S^{(\pm)})$ denote the two
Hopf algebras mentioned above. Then from Proposition \ref{deau}, we
conclude that $\A_q(n)$ is a left $D_q^{(\pm)}$-module algebra,
which gives rise to smash product algebras $\A_q(n)\sharp
D_q^{(\pm)}$ (see \cite{Hu1}).
\begin{defn}$($\cite[Def. 3.5]{Hu1}$)$
Let $\mathcal{W}_{q}(2n)$ be the associative algebra over~$\bK$
generated by the symbols
$\T(\pm\ve_i),\si_i^{\pm1},x_{i},\partial_{i} ~(1\le i\le n)$ with
the following defining relations:
\begin{eqnarray*}
&&\T(\pm\ve_i)\T(\mp\ve_i)=1=\si_i^{\pm1}\si_i^{\mp1},\quad
\Theta({-}\ve_i{+}\ve_{i+1})=\si_i \si_{i+1},\\
&&\Theta(\vei)\Theta(\ve_j)=\Theta(\vei{+}\ve_j)=\Theta(\ve_j)\Theta(\vei),\quad
\si_i\si_j=\si_j\si_i,\\
&&\si_i\Theta(\ve_j)=\Theta(\ve_j)\si_i,\quad
\Theta(\vei)x_{j}\Theta({-}\vei)=\theta(\vei,\ve_j)x_{j},\\
&&\Theta(\vei)\p_{j}\Theta({-}\vei)=\theta(\ve_j,\vei)\p_{j},\quad
\si_i x_{j}\si_i^{-1}=q^{\delta_{ij}}x_{j},\\
&&\si_i\p_{j}\si_i^{-1}=q^{-\delta_{ij}}\p_{j},\quad
x_{i}x_{j}=\theta(\vei,\ve_j)x_{j}x_{i},\\
&&\p_i\p_j=\theta(\vei,\ve_j)\p_j\p_i, \quad
\p_ix_{j}=\theta(\ve_j,\vei)x_{j}\p_i \quad (i\not=j),\quad
\p_ix_{i}-q^{\pm1}x_{i}\p_i=\si_i^{\mp1},
\end{eqnarray*}
where the last relations are equivalent to the relations:
$$\p_ix_{i}=\frac{q\si_i-(q\si_i)^{-1}}{q-q^{-1}},\qquad\qquad
x_{i}\p_i=\frac{\si_i-\si_i^{-1}}{q-q^{-1}}.$$
\end{defn}

\section{Realization of the simple root vectors}
\subsection{}
Based on the results above, we will realize the higher rank
quantized enveloping algebra $U_q(sl_{n+1})$ as some quantum
differential operators in $\W_q(2n)$ defined over the quantum
divided algebra $\A_q(n)$ such that $\A_q(n)$ is made into a
nonhomogeneous $U_q(sl_{n+1})$-module. This amounts to saying that
we can quantize the components of degree $-1$ and degree $0$ of
$\W(n)$, together with a $n$-dimensional subspace $L_1$ of
$\W(n)_1$.

%As we mentioned in the Introduction,
%it is worth pointing out that we realize the higher rank $U_q(sl_{n+1})$
%as certain quantum differential operators in $n$ variables rather than $n+1$
%variables.

First of all, we give the following lemma.
\begin{lem}\label{theta}
$\sum\limits_{k=1}^n\theta(\ve_k,\beta)[\beta_k]
=\sum\limits_{k=1}^n\theta\left(\ve_k,\beta{+}m(\vei{-}\ve_{i+1})
\right) [\left(\beta{+}m(\vei{-}\ve_{i+1})\right)_k], \ \forall\;m\in
\bZ$.
\end{lem}
\begin{proof}
 It suffices to verify
 \begin{eqnarray*}
\theta(\vei,\beta)[\bi]+\theta(\ve_{i+1},\beta)[\beta_{i+1}]
&=&\theta\left(\ve_i,\beta{+}m(\vei{-}\ve_{i+1})\right)[\left(\beta{+}m(\vei{-}\ve_{i+1})_i\right)]\\
&&+\theta\left(\ve_{i+1},\beta{+}m(\vei{-}\ve_{i+1})\right)[\left(\beta{+}m(\vei{-}\ve_{i+1})_{i+1}\right)].
\end{eqnarray*}
\begin{eqnarray*}
RHS&=&\theta(\vei,\beta)\theta(\vei,\vei)^m\theta(\vei,{-}\ve_{i+1})^m[\left(\beta{+}m(\vei{-}\ve_{i+1})_i\right)]\\
&&+\theta(\ve_{i+1},\beta)\theta(\ve_{i+1},\vei)^m\theta(\ve_{i+1},{-}\ve_{i+1})^m
[\left(\beta{+}m(\vei{-}\ve_{i+1})_{i+1}\right)]\\
&=&q^m\theta(\vei,\beta)[\bi{+}m]+q^m\theta(\ve_{i+1},\beta)[\beta_{i+1}{-}m]\\
&=&\theta(\vei,\beta)[\bi]+\theta(\ve_{i+1},\beta)[\beta_{i+1}]
+\left(\theta(\vei,\beta)q^{\bi}-\theta(\ve_{i+1},\beta)q^{-\beta_{i+1}}\right)q^m[m]\\
&=&LHS.
\end{eqnarray*}
This confirms the identity.
\end{proof}

\subsection{}
We are now in a position to state and prove one of our main results.
\begin{theorem}\label{qsl_n+1}
 For any monomial $x^{(\beta)}\in \mathcal{A}_q(n)$ and $1\le i \le n$, set
\begin{eqnarray*}
e_i(x^{(\beta )})&=&x_i{\partial}_{i+1}{\sigma}_i(x^{(\beta)}),\\
f_i(x^{(\beta)})&=&{{\sigma}_i}^{-1}
x_{i+1}{\partial}_{i}(x^{(\beta)}),\\
K_i(x^{(\beta)})&=&{\sigma}_i{\sigma}_{i+1}^{-1}(x^{(\beta)}),\\
e_n(x^{(\beta)})&=&\left(\prod\limits_{i=1}^{n-1}{\sigma}_i^{-1}\right)
\left(x_n\sum\limits_{i=1}^n x_i \partial_i
\Theta(\varepsilon_i)\right)\left(x^{(\beta )}\right),\\
f_n(x^{(\beta)})&=&\left(-\partial_n\prod\limits_{i=1}^{n-1}\sigma_i\right)\left(x^{(\beta
)}\right),\\
K_n(x^{(\beta)})&=&\left(\sigma_n\prod\limits_{i=1}^{n}\sigma_i\right)\left(x^{(\beta
)}\right).
\end{eqnarray*}
These formulas make $\mathcal{A}_q(n)$ into a nonhomogeneous
$U_q(sl_{n+1})$-module.
\end{theorem}

\begin{rem}
 Taking $q=1$,
 we recover the realization of $sl_{n+1}$ (cf. Lemma \ref{sln+1}).
\end{rem}

\begin{proof}
By Proposition \ref{qsln}, we only need to check the algebraic
relations that the generators $\{\,e_n,\,f_n,\,K_n^{\pm1}\,\}$
occur.

Using Lemma \ref{*product and theta} and formulas (\ref{sigma}),
(\ref{partial}) and (\ref{theta}), we obtain
\begin{eqnarray*}
e_i(x^{(\beta )})&=& q^{\bi-\ve_{i+1}*\beta+\vei*(\beta-\ve_{i+1})}
\begin{bmatrix}
\beta{+}\vei{-}\ve_{i+1}\\ \vei
\end{bmatrix}
x^{(\beta+\vei-\ve_{i+1})}\\
&=&[\bi{+}1]x^{(\beta+\vei-\ve_{i+1})},\\
f_i(x^{(\beta)})&=&q^{-\vei*\beta+\ve_{i+1}*(\beta-\vei)-(\bi-1)}
\begin{bmatrix}
\beta{-}\vei{+}\ve_{i+1}\\ \ve_{i+1}
\end{bmatrix}
x^{(\beta-\vei+\ve_{i+1})}\\
&=&[\beta_{i+1}{+}1]x^{(\beta-\vei+\ve_{i+1})},\\
e_n(x^{(\beta)})&=&[\beta_n{+}1]\left(\sum\limits_{i=1}^n\theta(\vei,\beta)[\bi]\right)x^{(\beta+\ve_n)},\\
f_n(x^{(\beta)})&=&-x^{(\beta-\ve_n)},\quad K_n(x^{(\beta)})=q^{\sum\limits_{i=1}^n\bi+\beta_n}
x^{(\beta)}.
\end{eqnarray*}

First of all, it is clear that (R1) holds.

For (R2), we have
\begin{eqnarray*}
K_i e_n{K_i}^{-1}(x^{(\beta )})&=&q^{\beta_{i+1}-\bi}[\beta_{n+1}]
\left(\sum\limits_{i=1}^n\theta(\vei,\beta)[\bi]\right)K_ix^{(\beta+\ve_n)}\\
&=&q^{\beta_{i+1}-\bi}
q^{\bi-\beta_{i+1}+\delta_{i,n}-\delta_{i+1,n}}e_n(x^{(\beta)})\\
&=&q^{\delta_{i,n}-\delta e_{i+1,n}}e_n(x^{(\beta)})\\
&=&q^{a_{i,n}}e_n(x^{(\beta)}),
\end{eqnarray*}
for $1\le i \le n{-}1$.
\begin{eqnarray*}
K_n e_n{K_n}^{-1}(x^{(\beta)})
&=&q^{-(\sum\limits_{i=1}^n\bi+\beta_n)}K_n e_n(x^{(\beta)})\\
&=&q^{-(\sum\limits_{i=1}^n\bi+\beta_n)}[\beta_{n+1}]
\left(\sum\limits_{i=1}^n\theta(\vei,\beta)[\bi]\right)K_n
x^{(\beta+\ve_n)}\\
&=&q^{-(\sum\limits_{i=1}^n\bi+\beta_n)}
q^{(\sum\limits_{i=1}^n\bi+\beta_n)+2}e_n(x^{(\beta)})\\
&=&q^{a_{n,n}}e_n(x^{(\beta)}).
\end{eqnarray*}
Similarly, for $1\le i\le n{-}1$, we obtain
\begin{eqnarray*}
K_i f_n K_i^{-1}(x^{(\beta)})
&=&-q^{\beta_{i+1}-\bi}q^{-\beta_{i+1}+\bi+\delta_{i+1,n}}x^{(\beta-\ve_n)}\\
&=&q^{\delta_{i+1,n}}f_n(x^{(\beta)})
=q^{-a_{i,n}}f_n(x^{(\beta)}).
\end{eqnarray*}
For $i=n$, we have
\begin{eqnarray*}
K_n f_n{K_n}^{-1}(x^{(\beta)})
=-q^{-(\sum\limits_{i=1}^n\bi+\beta_n)}q^{(\sum\limits_{i=1}^n\bi+\beta_n)-2}(x^{(\beta)})=q^{-a_{nn}}f_n(x^{(\beta)}).
\end{eqnarray*}

As for (R3), firstly, let us prove $[e_n,f_i]=0$, for $1\le i\le n{-}1$. In fact, we have
\begin{eqnarray*}
[e_n,f_i](x^{(\beta)})&=&e_nf_i(x^{(\beta)})-f_ie_n(x^{(\beta)})\\
&=&[\beta_{i{+}1}{+}1][\beta_n{+}\delta_{i{+}1,n}{+}1]\left(\sum\limits_{k=1}^n\theta(\ve_k,\beta{-}\vei{+}\ve_{i{+}1})
[(\beta{-}\vei{+}\ve_{i{+}1})_k]\right)x^{(\beta{-}\vei{+}\ve_{i+1}{+}\ve_n)}\\
&&-[\beta_n{+}1][\beta_{i+1}{+}\delta_{i{+}1,n}{+}1]
\left(\sum\limits_{i=1}^n\theta(\vei,\beta)
[\bi]\right)x^{(\beta{-}\vei{+}\ve_{i{+}1}{+}\ve_n)}
=0,
\end{eqnarray*}
where the last equality is a consequence of Lemma \ref{theta}.

While for the relations $[e_i,f_n]=0,~ (1\le i\le n{-}1)$, these follows from
\begin{eqnarray*}
[e_i,f_n](x^{(\beta)})&=&e_if_n(x^{(\beta)})-f_ne_i(x^{(\beta)})\\
&=&-[\beta_i{+}1]x^{(\beta{+}\vei{-}\ve_{i{+}1}{-}\ve_n)}
-(-1)[\beta_i{+}1]x^{(\beta{+}\vei{-}\ve_{i{+}1}{-}\ve_n)}
=0.
\end{eqnarray*}

Finally, we shall check that $[e_n, f_n]=\fr{K_n-K_n^{-1}}{q-q^{-1}}$. Indeed,
\begin{eqnarray*}
[\,e_n,&f_n&](x^{(\beta)})=e_nf_n(x^{(\beta)})-f_ne_n(x^{(\beta)})\\
&=&-[\beta_n]
\left(\sum \limits_{i=1}^n\theta(\vei,\beta{-}\ve_n)[(\beta{-}\ve_n)_i]\right)x^{(\beta)}
-(-1)[\beta_n{+}1]
\left(\sum\limits_{i=1}^n\theta(\vei,\beta)[\beta_i]\right)x^{(\beta)}\\
&=&-\left(\sum\limits_{i=1}^{n-1}q\theta(\vei,\beta)[\beta_i][\beta_n]
+\theta(\ve_n,\beta)[\beta_n{-}1][\beta_n]\right)x^{(\beta)}\\
&&+\left(\sum\limits_{i=1}^{n-1}\theta(\vei,\beta)[\beta_i][\beta_n{+}1]
+\theta(\ve_n,\beta)[\beta_n][\beta_n{+}1] \right) x^{(\beta)}\\
&=&\left(\sum\limits_{i=1}^n
q^{-\beta_n}\theta(\vei,\beta)[\beta_i]+q^{\beta_n}\theta(\ve_n,\beta)[\beta_n]\right)x^{(\beta)}\\
&=&\fr{x^{(\beta)}}{q{-}q^{{-}1}}\left(\sum\limits_{i=1}^n
q^{{-}\beta_n}\theta(\vei,\beta)q^{\bi}{+}q^{\beta_n}\theta(\ve_n,\beta)q^{\beta_n}{-}
\sum\limits_{i=1}^n q^{{-}\beta_n}\theta(\vei,\beta)q^{{-}\bi}
{-}q^{\beta_n}\theta(\ve_n,\beta)q^{{-}\beta_n}\right)\\
&=&\fr{1}{q{-}q^{-1}}\left(q^{\sum\limits_{i=1}^n
\bi{+}\beta_n}-q^{{-}\sum\limits_{i=1}^n\bi{-}\beta_n}\right)x^{(\beta)}
=\fr{K_n-K_n^{-1}}{q-q^{-1}}x^{(\beta)}.
\end{eqnarray*}

As for (R4),
it is easily seen that
\begin{eqnarray*}
e_{n-1}^2e_n(x^{(\beta)})&=&[\beta_n{+}1][\beta_{n{-}1}{+}1][\beta_{n{-}1}{+}2]
\left(\sum\limits_{i=1}^n\theta(\vei,\beta)[\bi]\right)x^{(\beta+2\ve_{n-1}{-}\ve_n)},\\
e_{n-1}e_ne_{n-1}(x^{(\beta)})&=&[\beta_n][\beta_{n{-}1}{+}1][\beta_{n{-}1}{+}2]
\left(\sum\limits_{i=1}^n\theta(\vei,\beta)[\bi]\right)x^{(\beta{+}2\ve_{n-1}{-}\ve_n)},\\
e_ne_{n-1}^2(x^{(\beta)})&=&[\beta_n{-}1][\beta_{n-1}{+}1][\beta_{n-1}{+}2]
\left(\sum\limits_{i=1}^n\theta(\vei,\beta)[\bi]\right)x^{(\beta{+}2\ve_{n-1}{-}\ve_n)},
\end{eqnarray*}
where the last equality follows from Lemma 2.1 when taking $m=2$ and $i=n{-}1$. By definition,
$[m+1]-(q+q^{-1})[m]+[m-1]=0$,
we conclude that
$$e_{n-1}^2e_n-(q+q^{-1})e_{n-1}e_ne_{n-1}+e_ne_{n-1}^2=0.$$
Similarly, we have that
$$e_n^2e_{n-1}-(q+q^{-1})e_ne_{n-1}e_n+e_{n-1}e_n^2=0.$$

We are now in a position to show (R5),
namely, $e_ie_n=e_ne_i$ for $1\le i\le n{-}2$.
\begin{eqnarray*}
e_ie_n(x^{(\beta)})
&=&[\beta_n{+}1]\left(\sum\limits_{i=1}^n\theta(\vei,\beta)
[\bi]\right)e_ix^{(\beta{+}\ve_n)}\\
&=&[\beta_n{+}1][\bi{+}1]\left(\sum\limits_{i=1}^n\theta(\vei,\beta)
[\bi]\right)x^{(\beta{+}\ve_n{+}\vei{-}\ve_{i{+}1})}\\
&=&[\beta_n{+}1][\bi{+}1]\left(\sum\limits_{k=1}^n
\theta(\ve_k,\beta{+}\vei{-}\ve_{i{+}1})
[(\beta{+}\vei{-}\ve_{i{+}1})_k]\right)x^{(\beta{+}\ve_n{+}\vei{-}\ve_{i{+}1})}\\
&=&[\bi{+}1]e_n(x^{(\beta{+}\vei{-}\ve_{i{+}1})})\\ &=&e_ne_i(x^{(\beta)}),
\end{eqnarray*}
where the third equality is a consequence of Lemma \ref{theta}.

We will prove the relation (R6) as follows.
Firstly, it is not difficult to check that
\begin{eqnarray*}
f_{n-1}^2f_n(x^{(\beta)})&=&-[\beta_n][\beta_n{+}1]x^{(\beta{-}2\ve_{n{-}1}{+}\ve_n)},\\
f_{n-1}f_nf_{n-1}(x^{(\beta)})&=&-[\beta_n{+}1][\beta_n{+}1]x^{(\beta{-}2\ve_{n{-}1}{+}\ve_n)},\\
f_nf_{n-1}^2(x^{(\beta)})&=&-[\beta_n{+}2][\beta_n{+}1]x^{(\beta{-}2\ve_{n{-}1}{+}\ve_n)}.
\end{eqnarray*}
It is easy to see that
$f_{n-1}^2f_n-(q+q^{-1})f_{n-1}f_nf_{n-1}+f_nf_{n-1}^2=0$.
Likewise, we can show that
$f_n^2f_{n-1}-(q+q^{-1})f_nf_{n-1}f_n+f_{n-1}f_n^2=0$.

Finally, we check the relation (R7) as follows.
$$
f_if_n(x^{(\beta)})
=-[\beta_{i+1}{+}1]x^{(\beta{-}\ve_n{-}\vei{+}\ve_{i+1})}
=f_nf_i(x^{(\beta)}),
$$
for $1\le i\le n-2$.

This completes the proof.
\end{proof}

\section{Realization of the non-simple root vectors}
In this section,
based on the realization model of the simple root vectors,
we will use Lusztig symmetries to give the quantum differential operators realization
for the non-simple quantum root vectors of $U_q(sl_{n+1})$.
One will see that the expression of quantum root vectors by
our quantum differential operators is well-matched with the action of Lusztig symmetries. This implies further that
our construction is intrinsically natural.

\subsection{}
Let $\alpha_i=\varepsilon_i-\varepsilon_{i+1}(1\leq i\leq n-1)$ be
the simple roots of $sl_n$ and $s_i~(1\leq i\leq n-1)$ be the
corresponding reflections of the Weyl group $W=S_n$. Let $T_i$ be
the braid automorphism of $U_q(sl_n)$ determined by $s_i$ introduced
as $T_{i,-1}^{''}$ by Lusztig in \cite[\S37.1.3]{Lus}. They are
defined as follows:
\begin{eqnarray}
T_i(K_\mu)&=&K_{s_i(\mu)},\quad T_i(e_i)=-f_iK_i^{-1},\quad T_i(f_i)=-K_ie_i,\label{Tiei}\\
T_i(e_j)&=&e_j,\quad T_i(f_j)=f_j, \quad (|i-j|>1), \\
T_i(e_j)&=&e_ie_j-qe_je_i,\quad T_i(f_j)=f_jf_i-q^{-1}f_if_j,
\quad(|i-j|=1),\label{Tiei+1}
 \end{eqnarray}
where $e_i$, $f_i$ are the simple root vectors of $U_q(sl_n)$
associated to $\alpha_i$ and $-\alpha_i$.

Take a reduced presentation of the longest element
$w_0=s_{i_1}s_{i_2}\cdots s_{i_m}$ in $W$, then
\begin{eqnarray}
 e_\alpha=T_{i_1}\cdots T_{i_{p-1}}(e_{i_p}),
 \quad f_\alpha=T_{i_1}\cdots T_{i_{p-1}}(f_{i_p})
\end{eqnarray}
are called quantum root vectors of~ $U_q(sl_n)$ corresponding to
$\pm\alpha=\pm s_{i_1}s_{i_2}\cdots
s_{i_{p-1}}(\alpha_{i_p}),~ 1\le p \le m$.
From now on, let us fix a reduced presentation of $w_0$
as follows
$$w_0=s_1s_2s_1s_3s_2s_1\dots
s_{n-1}s_{n-2}\cdots s_2s_1$$
which gives a convex ordering of the positive roots system of
$sl_n$
$$\alpha_{12}, \,\alpha_{13}, \,\alpha_{23},\,
\alpha_{14},\, \alpha_{24},\, \alpha_{34},\,
\cdots,\, \alpha_{1n},\, \alpha_{2n},\, \cdots,\, \alpha_{n-1,n}.$$
Hence, all the quantum root vectors of~$U_q(sl_n)$ associated to this ordering
are determined. For instance,
\begin{eqnarray}
 e_{\alpha_{1n}}=T_1T_2T_1T_3T_2T_1\cdots T_{n-2}\cdots
T_2T_1(e_{n-1}).
\end{eqnarray}
Later on, we will use the quantum root vectors of
 $U_q(sl_n)$ associated to this ordering.

\begin{lem}$($\cite[\S39.2.4]{Lus}$)$\label{T_i}
The braid automorphisms
$T_i$'s satisfy the relations
 \begin{eqnarray*}
T_iT_jT_i&=&T_jT_iT_j,\quad |i-j|=1,\\
 T_iT_j(e_i)&=&e_j,\quad |i-j|=1,\\
T_i(e_j)&=&e_j,\quad |i-j|>1.
  \end{eqnarray*}
\end{lem}

\subsection{}
Now let us introduce some new $q$-differential operators in $\W_q(2n)$.

Set
\begin{eqnarray}
e_{s,n+1}&=&\left(\prod\limits_{i\neq s}{\sigma}_i^{-1}\right)
\left(x_s\sum\limits_{i=1}^n x_i \partial_i
\Theta(\varepsilon_i)\right), ~ 1\le s \le n,\label{e_sn+1}\\
e_{n+1,s}&=&-\partial_s\left(\prod\limits_{i\neq s}{\sigma}_i\right),~ 1\le s\le n.\label{e_n+1s}
\end{eqnarray}

Then we have
\begin{eqnarray}
 e_{s,n+1}(x^{(\beta)})
 &=&\left(\prod\limits_{i\neq
s}{\sigma}_i^{-1}\right) \left(x_s\sum\limits_{i=1}^n x_i \partial_i
\Theta(\varepsilon_i)\right)(x^{(\beta)})\label{esn+1}\\
&=&\left(\prod\limits_{i\neq
s}{\sigma}_i^{-1}\right)x_s\left(\sum\limits_{i=1}^n
\theta(\varepsilon_i,\beta)[\beta_i]\right)(x^{(\beta)}) \nonumber\\
&=&q^{-\sum\limits_{i=s+1}^n\beta_i}[\beta_s{+}1]\left(\sum\limits_{i=1}^n
\theta(\varepsilon_i,\beta)[\beta_i]\right)(x^{(\beta+\varepsilon_s)}), \nonumber\\
e_{n+1,s}(x^{(\beta)})&=&-\partial_s\left(\prod\limits_{i\neq s}{\sigma}_i\right)(x^{(\beta)})
=-q^{\sum\limits_{i=s+1}^n\beta_i}(x^{(\beta-\varepsilon_s)})\label{en+1s}.
\end{eqnarray}

In the remainder of this section, we will show that
the $q$-differential operators defined above are coincident with the one of four kinds of
Lusztig's quantum root vectors.

It is not difficult to check that
\begin{eqnarray}
 e_{i,j}(x^{(\beta)})
 &=&q^{-\sum\limits_{i<s<j}\beta_s}[\beta_i{+}1]
 x^{(\beta{+}\varepsilon_i{-}\varepsilon_j)},\quad(1\leq i<j\leq n)\label{eij},\\
e_{i,j}(x^{(\beta)})&=&q^{\sum\limits_{j<s<i}\beta_s}
[\beta_i{+}1]x^{(\beta{-}\varepsilon_j{+}\varepsilon_i)},\quad \ \ \,(1\leq j<i\leq n)\label{eji}.
\end{eqnarray}

\begin{prop}\label{qbracket}
For any $1\leq s<n,s<j\leq n$,
one has
\begin{eqnarray*}
&&(1)~~ e_{s,n+1}=[e_{s,s+1},e_{s+1,n+1}]_q:=e_{s,s+1}e_{s+1,n+1}-qe_{s+1,n+1}e_{s,s+1},\\
 &&(2)~~ e_{s,n+1}=[e_{s,j},e_{j,n+1}]_q:=e_{s,j}e_{j,n+1}-qe_{j,n+1}e_{s,j},\\
 &&(3)~~ e_{n+1,s}=[e_{n+1,j},e_{j,s}]_{q^{-1}},\\
 &&(4)~~ [e_{s,n+1},e_{n+1,s}]=\frac{\left(\prod_{i=1}^{n} \sigma_i\right)\sigma_s-\left(\prod_{i=1}^{n}
 \sigma_i^{-1}\right)\sigma_s^{-1}}{q-q^{-1}}, \quad \mbox{where}~~\left(\prod_{i=1}^{n}
 \sigma_i\right)\sigma_s=K_{\varepsilon_s-\varepsilon_{n+1}}.
  \end{eqnarray*}
\end{prop}
This proposition asserts that the expressions of
 $e_{s,n+1}$ and $e_{n+1,s}$ by the $q$-brackets
 are independent of the choice of $j~ (s<j\le n)$.

\begin{proof}
(1) \ For $1\leq s<n$, using formulas (\ref{esn+1}) and (\ref{eij}), we have
\begin{eqnarray*}
\left(e_{s,s+1}e_{s+1,n+1}\right.&-&\left.q\,e_{s+1,n+1}e_{s,s+1}\right)(x^{(\beta)})\\
\quad &=&q^{-\sum\limits_{i=s+2}^n\beta_i}[\beta_{s+1}{+}1]\left(\sum\limits_{i=1}^n
\theta(\varepsilon_i,\beta)[\beta_i]\right)e_{s,s+1}(x^{(\beta+\varepsilon_{s+1})})\\
\quad&&-q\,[\beta_s{+}1]e_{s+1,n+1}(x^{(\beta+\varepsilon_{s}-\varepsilon_{s+1})})\\
\quad&=&q^{-\sum\limits_{i=s+2}^n\beta_i}[\beta_{s+1}{+}1]\left(\sum\limits_{i=1}^n
\theta(\varepsilon_i,\beta)[\beta_i]\right)[\beta_s{+}1](x^{(\beta+\varepsilon_s)})\\
\quad&&-q\,[\beta_s{+}1]q^{-\sum\limits_{i=s+2}^n\beta_i}[\beta_{s+1}]\left(\sum\limits_{i=1}^n
\theta(\varepsilon_i,\beta{+}\varepsilon_s{-}\varepsilon_{s+1})[(\beta{+}\varepsilon_s{-}\varepsilon_{s+1})_i]\right)
(x^{(\beta+\varepsilon_s)})\\
\quad&=&q^{-\sum\limits_{i=s+1}^n\beta_i}[\beta_s{+}1]\left(\sum\limits_{i=1}^n
\theta(\varepsilon_i,\beta)[\beta_i]\right)(x^{(\beta+\varepsilon_s)})\\
\quad&=&e_{s,n+1}(x^{(\beta)}).
\end{eqnarray*}

(2) \ We will prove (2) based on (1) and by induction.

From (1), we know that the formula holds for $s=n{-}1$.

Assume that $e_{j+1,n+1}=[e_{j+1,k},e_{k,n+1}]_q$ hold
for any $j{+}1$ with $j{+}1<k<n{+}1$. We will prove it  for $j$ as follows.
\begin{eqnarray*}
e_{j,n+1}&=&[e_{j,j+1},e_{j+1,n+1}]_q=e_{j,j+1}e_{j+1,n+1}-qe_{j+1,n+1}e_{j,j+1}\\
&=&e_{j,j+1}(e_{j+1,k}e_{k,n+1}-qe_{k,n+1}e_{j+1,k})-q(e_{j+1,k}e_{k,n+1}-qe_{k,n+1}e_{j+1,k})e_{j,j+1}\\
&=&(e_{j,j+1}e_{j+1,k}-qe_{j+1,k}e_{j,j+1})e_{k,n+1}-qe_{k,n+1}(e_{j,j+1}e_{j+1,k}-qe_{j+1,k}e_{j,j+1})\\
&=&e_{j,k}e_{k,n+1}-qe_{k,n+1}e_{j,k}\\
&=&[e_{j,k},e_{k,n+1}]_q.
\end{eqnarray*}
(3)\quad For any ~$1\leq s<n,s<j\leq n$,
it follows from formulas (\ref{en+1s}) and (\ref{eji}) that
\begin{eqnarray*}
(e_{n+1,j}e_{j,s}&-&q^{-1}e_{j,s}e_{n+1,j})(x^{(\beta)})\\
&=&[\beta_j{+}1]q^{\sum\limits_{j>i>s}\beta_i}e_{n+1,j}(x^{(\beta-\varepsilon_s+\varepsilon_j)})
+q^{-1}\cdot q^{\sum\limits_{i=j+1}^n\beta_i}e_{j,s}(x^{(\beta-\varepsilon_j)})\\
&=&-[\beta_j{+}1]q^{\sum\limits_{j>i>s}\beta_i}q^{\sum\limits_{i=j+1}^n\beta_i}(x^{(\beta-\varepsilon_s)})
+ q^{\sum\limits_{i=j+1}^n\beta_i-1}[\beta_j]q^{\sum\limits_{j>i>s}\beta_i}(x^{(\beta-\varepsilon_s)})\\
&=&-q^{\sum\limits_{i=s+1}^n\beta_i}(x^{(\beta-\varepsilon_s)})=e_{n+1,s}(x^{(\beta)}).
\end{eqnarray*}
(4)\quad For any $1\leq s<n$, applying formulas (\ref{esn+1}) and (\ref{en+1s}),
we obtain
\begin{eqnarray*}
(e_{s,n+1}e_{n+1,s}&-&e_{n+1,s}e_{s,n+1})(x^{(\beta)})\\
&=&-q^{\sum\limits_{i=s+1}^n\beta_i}e_{s,n+1}(x^{(\beta-\varepsilon_s)})
-q^{-\sum\limits_{i=s+1}^n\beta_i}[\beta_s{+}1]\left(\sum\limits_{i=1}^n
\theta(\varepsilon_i,\beta)[\beta_i]\right)e_{n+1,s}(x^{(\beta+\varepsilon_s)})\\
&=&-[\beta_s]\left(\sum\limits_{i=1}^n
\theta(\varepsilon_i,\beta{-}\varepsilon_s)[(\beta{-}\varepsilon_s)_i]\right)(x^{(\beta)})
+[\beta_s{+}1]\left(\sum\limits_{i=1}^n
\theta(\varepsilon_i,\beta)[\beta_i]\right)(x^{(\beta)})\\
&=&-[\beta_s]\left(
\sum\limits_{i=1}^{s-1}
\theta(\varepsilon_i,\beta{-}\varepsilon_s) [\beta_i]+
\theta(\varepsilon_s,\beta)[\beta_s{-}1]
+\sum\limits_{i=s+1}^n
\theta(\varepsilon_i,\beta{-}\varepsilon_s)[\beta_i]
\right)(x^{(\beta)})\\
&&+[\beta_s{+}1]\left(
\sum\limits_{i=1}^{s-1}
\theta(\varepsilon_i,\beta) [\beta_i]+
\theta(\varepsilon_s,\beta) [\beta_s]
+\sum\limits_{i=s+1}^n
\theta(\varepsilon_i,\beta) [\beta_i]
\right)(x^{(\beta)})\\
&=&\sum\limits_{i=1}^{s-1}
\left(([\beta_s]q+q^{-\beta_s})\theta(\varepsilon_i,\beta) [\beta_i]-[\beta_s]q\theta(\varepsilon_i,\beta) [\beta_i]
\right)(x^{(\beta)})\\
&&+\left(([\beta_s]q^{-1}+q^{\beta_s})\theta(\varepsilon_s,\beta) [\beta_s]-[\beta_s]\theta(\varepsilon_s,\beta) ([\beta_s]q^{-1}-q^{-\beta_s})\right)(x^{(\beta)})\\
&&+\sum\limits_{i=s+1}^n\left(
([\beta_s]q^{-1}+q^{\beta_s})\theta(\varepsilon_i,\beta) [\beta_i]-[\beta_s]q^{-1}\theta(\varepsilon_i,\beta) [\beta_i]
\right)(x^{(\beta)})\\
&=&\left(q^{-\sum\limits_{i=s}^{n}\beta_i}\left[\,\sum\limits_{i=1}^{s}\beta_i\,\right]
+q^{\sum\limits_{i=1}^{s}\beta_i}\left[\,\sum\limits_{i=s}^{n}\beta_i\,\right]
\right)(x^{(\beta)})
=\left[\,\beta_s{+}\sum\limits_{i=1}^{n}\beta_i\,\right](x^{(\beta)})\\
&=&\frac{\left(\prod_{i=1}^{n} \sigma_i\right)\sigma_s-\left(\prod_{i=1}^{n}
 \sigma_i^{-1}\right)\sigma_s^{-1}}{q-q^{-1}}(x^{(\beta)}).
\end{eqnarray*}

Thus, we complete the proof.
\end{proof}

\subsection{}
We recall that $e_{ij}=x_i\p_j\si_i~(1\le i <j \le n)$,
$e_{ij}=\si_j^{-1}x_i\p_j~(1\le j <i\le n)$
as given in the Introduction and
$e_{s,n+1},e_{n+1,s}$ are defined by formulas (\ref{e_sn+1}) and (\ref{e_n+1s}).

\begin{theorem}\label{root vector}
If we identify  the generators of $U_q(sl_{n+1})$ with certain
$q$-differential operators in $\W_q(2n)$,
namely,
$e_i:=e_{i,i+1},\, f_i:=e_{i+1,i},\, K_i:=\sigma_i\sigma_{i+1}^{-1}, \ 1\le i<n$, and $e_n:=e_{n,n{+}1}$, $f_n:=e_{n{+}1,n}$, $K_n:=\sigma_n\prod_{i=1}^n\sigma_i$, then we have

(1)~~ $e_{\alpha_{ij}}$ correspond to $e_{i,j}$  ~$(1\leq i<j\leq n{+}1)$,
where $e_{\alpha_{ij}}$
are positive root vectors
associated to $\alpha_{ij}=\varepsilon_i-\varepsilon_j$.

(2)~~ $f_{\alpha_{ji}}$ correspond to $e_{i,j}$ ~$(1\leq j<i\leq n{+}1)$,
where $f_{\alpha_{ji}}$
are negative root vectors
associated to $\alpha_{ji}=\varepsilon_i-\varepsilon_j$.
\end{theorem}

In order to prove the theorem,
we first establish an auxiliary lemma.

\begin{lem}\label{auxlem}
 With the identification as in Theorem \ref{root vector}, one has
\begin{eqnarray*}
&(1)&~\mbox {If}~~2\leq s{+}1<j\leq n{+}1,~ \mbox {then}~ [e_{sj},T_s(e_s)]_q=e_{s+1,j}.\\
&(2)&~\mbox {If}~~2\leq s{+}1<j\leq n{+}1, ~\mbox {then}~[T_s(f_s),e_{js}]_{q^{-1}}=e_{j,s+1}.
  \end{eqnarray*}
\end{lem}
\begin{proof}
(1)\quad
From \cite[Lemma 4.6]{Hu1}, it is clear that the formula holds for $2\le s{+}1< j \le n$. Hence it suffices to show the case $j=n{+}1$.
Using formulas (\ref{Tiei}), (\ref{esn+1}) and (\ref{eji}), we have
\begin{eqnarray*}
 LHS(&x^{(\beta)}&)=(qf_s{K_s}^{-1}e_{s,n+1}-e_{s,n+1}f_s{K_s}^{-1})(x^{(\beta)})\\
 &=&q^{{-}\sum\limits_{i=s+1}^n\beta_i{+}1}[\beta_s{+}1]\left(\sum\limits_{i=1}^n
\theta(\varepsilon_i,\beta)[\beta_i]\right)f_s{K_s}^{-1}(x^{(\beta+\varepsilon_s)})\\
&&-q^{{-}\beta_s{+}\beta_{s+1}}[\beta_{s+1}{+}1]e_{s,n+1}
(x^{(\beta{-}\varepsilon_s{+}\varepsilon_{s+1})})\\
&=&q^{{-}\sum\limits_{i=s+1}^n\beta_i{+}1}[\beta_s{+}1]\left(\sum\limits_{i=1}^n
\theta(\varepsilon_i,\beta)[\beta_i]\right)
q^{{-}\beta_s{-}1{+}\beta_{s+1}}[\beta_{s+1}{+}1]
(x^{(\beta+\varepsilon_{s+1})})\\
&&-q^{-\beta_s+\beta_{s+1}}[\beta_{s+1}{+}1]
q^{{-}\sum\limits_{i=s+1}^n\beta_i{-}1}[\beta_s]
\left(\sum\limits_{i=1}^n
\theta(\varepsilon_i,\beta{-}\varepsilon_s{+}\varepsilon_{s+1})
[(\beta{-}\varepsilon_s{+}\varepsilon_{s+1})_i]\right)(x^{(\beta+\varepsilon_{s+1})})\\
&=&q^{{-}\sum\limits_{i=s+2}^n\beta_i}[\beta_{s+1}{+}1]\left(\sum\limits_{i=1}^n
\theta(\varepsilon_i,\beta)[\beta_i]\right)(x^{(\beta+\varepsilon_{s+1})})\\
&&+q^{{-}\beta_s{+}\beta_{s+1}}[\beta_{s+1}{+}1]q^{{-}\sum\limits_{i=s+1}^n\beta_i{-}1}[\beta_s]\\
&&\cdot\left((\sum\limits_{i=1}^n
\theta(\varepsilon_i,\beta)[\beta_i])-(\sum\limits_{i=1}^n
\theta(\varepsilon_i,\beta{-}\varepsilon_s{+}\varepsilon_{s+1})[(\beta{-}\varepsilon_s{+}\varepsilon_{s+1})_i])\right)
(x^{(\beta+\varepsilon_{s+1})})\\
&=&e_{s+1,n+1}(x^{(\beta)})=RHS(x^{(\beta)}).
 \end{eqnarray*}
 Similarly, one can prove (2), by
formulas (\ref{Tiei}), (\ref{en+1s}) and (\ref{eij}).
\end{proof}

\subsection{}
We proceed to show Theorem \ref{root vector}.
\begin{proof}
It suffices to show (1). (2) can be proved in the same manner.

From Proposition \ref{prop4.6}, we conclude that the assertion is true for
 $1\le i, \,j \le n$. Therefore, we need only to deal with the case
 $j=n{+}1$, namely, $e_{\alpha_{s,n+1}}:=e_{s,n+1}$.

For $k=1$, since $T_i$ are automorphisms of $U_q(sl_{n+1})$, it follows from
Proposition \ref{prop4.6}, Lemma \ref{T_i} and formula (\ref{Tiei+1}) that
\begin{eqnarray*}
e_{\alpha_{1,n+1}}&=&T_1T_2T_1T_3T_2T_1\cdots T_{n-1}\cdots
T_2T_1(e_n)\\
&=&T_1T_2T_1T_3T_2T_1\cdots T_{n-2}\cdots
T_2T_1T_{n-1}(e_n)\\
&=&T_1T_2T_1T_3T_2T_1\cdots T_{n-2}\cdots
T_2T_1([e_{n-1},e_n]_q)\\
&=&[T_1T_2T_1T_3T_2T_1\cdots T_{n-2}\cdots
T_2T_1(e_{n-1}),e_n]_q\\
&=&[e_{1,n},e_n]_q=e_{1,n+1}.
\end{eqnarray*}
Assume that the claim holds for $j$, that is to say,
$$e_{\alpha_{j,n+1}}=T_1T_2T_1T_3T_2T_1\cdots T_{n-1}\cdots
T_2T_1T_nT_{n-1}\cdots T_{n-j+2}(e_{n-j+1}):=e_{j,n+1}.$$
We will prove it for $j{+}1$.
 Using Lemma \ref{T_i} and Lemma \ref{auxlem}, we obtain
\begin{eqnarray*}
\begin{split}
&e_{\alpha_{j+1,n+1}}
=T_1T_2T_1T_3T_2T_1\cdots T_{n-1}\cdots
T_2T_1T_nT_{n-1}\cdots
T_{n-j+2}T_{n-j+1}(e_{n-j})\\
&\ =T_1T_2T_1T_3T_2T_1\cdots T_{n-1}\cdots
T_2T_1T_nT_{n-1}\cdots
T_{n-j+2}([e_{n-j+1},e_{n-j}]_q)\\
&\ =[T_1T_2T_1T_3T_2T_1\cdots T_{n-1}\cdots T_2T_1T_nT_{n-1}\cdots
T_{n-j+2}(e_{n-j+1}),\\
&\qquad T_1T_2T_1T_3T_2T_1\cdots T_{n-1}\cdots
T_2T_1T_nT_{n-1}\cdots
T_{n-j+2}(e_{n-j})]_q\\
&\ =[e_{j,n+1},T_1T_2T_1T_3T_2T_1\cdots T_{n-1}T_{n-2}\cdots
T_{n-j-1}(e_{n-j})]_q\\
&\ =[e_{j,n+1},T_1T_2T_1T_3T_2T_1\cdots
T_{n-2}T_{n-3}\cdots
T_{n-j-2}(e_{n-j-1})]_q\\
&\quad=\cdots\cdots\\
&\ =[e_{j,n+1},T_1T_2T_1T_3T_2T_1\cdots T_{j+1}T_{j}\cdots
T_{2}T_1(e_2)]_q\\
&\ =[e_{j,n+1},T_1T_2T_1T_3T_2T_1\cdots T_{j}T_{j-1}\cdots
T_1(e_1)]_q\\
&\ =[e_{j,n+1},T_1T_2T_1T_3T_2T_1\cdots T_{j}T_{j-1}\cdots
T_2(-f_1 K_1^{-1})]_q\\
&\ =[e_{j,n+1},T_1T_2T_1T_3T_2T_1\cdots T_{j}T_{j-1}\cdots
T_2(-f_1)\cdot T_1T_2T_1T_3T_2T_1\cdots T_{j}T_{j-1}\cdots
T_2(K_1^{-1})]_q\\
&\ =[e_{j,n+1},-f_j K_j^{-1}]_q
=[e_{j,n+1},T_j(e_j)]_q
=e_{j+1,n+1}.
\end{split}
\end{eqnarray*}

This completes the proof.
\end{proof}

%\clearpage
\end{document}